# Error Estimation of Difference Operators on Irregular Nodes

H. Isshiki, Institute of Mathematical Analysis, Osaka, Japan,
isshiki@dab.hi-ho.ne.jp

T. Kawamura, Computational Fluid Dynamics Consulting Inc.,
kawamura@cfd-consulting.co.jp

D. Kitazawa, Institute of Industrial Science,
The University of Tokyo, dkita@iis.u-tokyo.ac.jp

**Abstract**
Error estimation of difference operators on irregular nodes is discussed. We can obtain the similar estimates of the errors. However, the error estimate for the difference operators for the second derivatives becomes lower because of asymmetric allocation of the nodal points.

**1. Imtroduction**

Moving Particle Semi-implicit method (MPS)[1] is a particle method and is widely used. It gives plausible numerical results in many cases. We have already given the mathematical base of difference operators in MPS[2]. Iribe and Nakaza has also proposed a method to improve the accuracy of the gradient operator[3].

Iribe-Nakaza's method is a special case of the author's method called Discrete Differential operators on Irregular Nodes (DDIN)[4]. In the author's paper[4], interpolation of discrete data given on irregular mesh is used. If we use the interpolation using power functions, the author's method can give a similar result obtained by Iribe and Nakaza for the gradient operator.

In the present paper, we discuss the error estimation of difference operators on irregular nodes theoretically. We also give the validation of the theory through numerical verification.

**2. Summary of error estimation on regular nodes**

In case of regular nodes, we can obtain the error estimates of the difference operators for the derivatives using Taylor expansion. The accuracy of the difference operators increases because of the symmetry of the allocation of the nodal points. The error estimates of the differential operators in regular nodes are summarized for 2D case as follows.

Central difference:
$$\left(\frac{\partial u}{\partial x}\right)_{i,j} = \frac{u_{i+1,j} - u_{i-1,j}}{2\Delta x} + O(\Delta x^2). \qquad (1)$$

Forward difference:
$$\left(\frac{\partial u}{\partial x}\right)_{i,j} = \frac{u_{i+1,j} - u_{i,j}}{2\Delta x} + O(\Delta x). \qquad (2)$$

Backward difference:
$$\left(\frac{\partial u}{\partial x}\right)_{i,j} = \frac{u_{i,j} - u_{i-1,j}}{2\Delta x} + O(\Delta x). \qquad (3)$$

Difference for the second derivatives:
$$\left(\frac{\partial^2 u}{\partial x^2}\right)_{i,j} = \frac{u_{i,j} - 2u_{i,j} + u_{i-1,j}}{\Delta x^2} + O(\Delta x^2). \qquad (4)$$

The formulas for $(\partial u/\partial y)_{i,j}$ and $(\partial^2 u/\partial y^2)_{i,j}$ could be obtained easily.

**3. Error estimation on regular nodes**

For simplicity, let's consider in 1D. The second order approximation of $\phi$ is given by

$$\phi_j = \phi_i + (\phi_x)_i(x_j - x_i) + \frac{1}{2}(\phi_{xx})_i(x_j - x_i)^2 + O(\varepsilon^3), \quad j = 1,2,\cdots J, \tag{5}$$

where

$$|x_j - x_i| < \varepsilon, \quad j = 1,2,\cdots,J. \tag{6}$$

Then, we have

$$\begin{bmatrix} x_1 - x_i & (x_1 - x_i)^2 \\ x_2 - x_i & (x_2 - x_i)^2 \\ \vdots & \vdots \\ \vdots & \vdots \\ x_J - x_i & (x_J - x_i)^2 \end{bmatrix} \begin{bmatrix} (\phi_x)_i \\ (\frac{1}{2}\phi_{xx})_i \end{bmatrix} = \begin{bmatrix} \phi_1 - \phi_i \\ \phi_2 - \phi_i \\ \vdots \\ \vdots \\ \phi_J - \phi_i \end{bmatrix} + \begin{bmatrix} O(\varepsilon^3) \\ O(\varepsilon^3) \\ \vdots \\ \vdots \\ O(\varepsilon^3) \end{bmatrix}. \tag{7}$$

Applying Least Square Method (LSM), we obtain

$$\begin{bmatrix} x_1 - x_i & x_2 - x_i & \cdots & \cdots & x_J - x_i \\ (x_1 - x_i)^2 & (x_2 - x_i)^2 & \cdots & \cdots & (x_J - x_i)^2 \end{bmatrix} \begin{bmatrix} x_1 - x_i & (x_1 - x_i)^2 \\ x_2 - x_i & (x_2 - x_i)^2 \\ \vdots & \vdots \\ \vdots & \vdots \\ x_J - x_i & (x_J - x_i)^2 \end{bmatrix} \begin{bmatrix} (\phi_x)_i \\ (\frac{1}{2}\phi_{xx})_i \end{bmatrix}$$

$$= \begin{bmatrix} x_1 - x_i & x_2 - x_i & \cdots & \cdots & x_J - x_i \\ (x_1 - x_i)^2 & (x_2 - x_i)^2 & \cdots & \cdots & (x_J - x_i)^2 \end{bmatrix} \begin{bmatrix} \phi_1 - \phi_i \\ \phi_2 - \phi_i \\ \vdots \\ \vdots \\ \phi_J - \phi_i \end{bmatrix} + \begin{bmatrix} O(\varepsilon^4) \\ O(\varepsilon^5) \end{bmatrix}. \tag{8}$$

Estimating the orders, we have

$$\begin{bmatrix} x_1 - x_i & x_2 - x_i & \cdots & \cdots & x_J - x_i \\ (x_1 - x_i)^2 & (x_2 - x_i)^2 & \cdots & \cdots & (x_J - x_i)^2 \end{bmatrix} \begin{bmatrix} x_1 - x_i & (x_1 - x_i)^2 \\ x_2 - x_i & (x_2 - x_i)^2 \\ \vdots & \vdots \\ \vdots & \vdots \\ x_J - x_i & (x_J - x_i)^2 \end{bmatrix}$$

$$= \begin{bmatrix} \sum (x_k - x_i)^2 & \sum (x_k - x_i)^3 \\ \sum (x_k - x_i)^3 & \sum (x_k - x_i)^4 \end{bmatrix} = \begin{bmatrix} O(\varepsilon^2) & O(\varepsilon^3) \\ O(\varepsilon^3) & O(\varepsilon^4) \end{bmatrix} \tag{9}$$

and

$$\begin{bmatrix} \sum (x_k - x_i)^2 & \sum (x_k - x_i)^3 \\ \sum (x_k - x_i)^3 & \sum (x_k - x_i)^4 \end{bmatrix}^{-1}$$

$$= \frac{1}{\sum (x_k - x_i)^2 \sum (x_k - x_i)^4 - \sum (x_k - x_i)^3 \sum (x_k - x_i)^3} \begin{bmatrix} \sum (x_k - x_i)^4 & -\sum (x_k - x_i)^3 \\ -\sum (x_k - x_i)^3 & \sum (x_k - x_i)^2 \end{bmatrix} \tag{10}$$

$$= \begin{bmatrix} O(\varepsilon^{-2}) & O(\varepsilon^{-3}) \\ O(\varepsilon^{-3}) & O(\varepsilon^{-4}) \end{bmatrix}.$$

Finally, we derive

$$
\begin{bmatrix} (\phi_x)_i \\ (\tfrac{1}{2}\phi_{xx})_i \end{bmatrix} = \left( \begin{bmatrix} x_1-x_i & x_2-x_i & \cdots & \cdots & x_J-x_i \\ (x_1-x_i)^2 & (x_2-x_i)^2 & \cdots & \cdots & (x_J-x_i)^2 \end{bmatrix} \begin{bmatrix} x_1-x_i & (x_1-x_i)^2 \\ x_2-x_i & (x_2-x_i)^2 \\ \vdots & \vdots \\ \vdots & \vdots \\ x_J-x_i & (x_J-x_i)^2 \end{bmatrix} \right)^{-1}
$$

$$
\cdot \begin{bmatrix} x_1-x_i & x_2-x_i & \cdots & \cdots & x_J-x_i \\ (x_1-x_i)^2 & (x_2-x_i)^2 & \cdots & \cdots & (x_J-x_i)^2 \end{bmatrix} \begin{bmatrix} \phi_1-\phi_i \\ \phi_2-\phi_i \\ \vdots \\ \vdots \\ \phi_J-\phi_i \end{bmatrix} + \begin{bmatrix} O(\varepsilon^2) \\ O(\varepsilon^1) \end{bmatrix}
\qquad (11)
$$

since

$$
\begin{bmatrix} O(\varepsilon^{-2}) & O(\varepsilon^{-3}) \\ O(\varepsilon^{-3}) & O(\varepsilon^{-4}) \end{bmatrix} \begin{bmatrix} O(\varepsilon^4) \\ O(\varepsilon^5) \end{bmatrix} = \begin{bmatrix} O(\varepsilon^2) \\ O(\varepsilon^1) \end{bmatrix}
\qquad (12)
$$

If we compare Eqs. (1) and (11), the accuracies of the first order difference operators are same unexpectedly. If we compare Eqs. (4) and (11), the accuracy becomes lower in irregular nodes than in regular nodes. The symmetry of the node allocation in regular nodes increases the accuracy.

## 4. Numerical examples

In 2D, the second order approximation of $\phi$ is given by

$$
\begin{bmatrix} u_1 & v_1 & u_1^2 & u_1 v_1 & v_1^2 \\ u_2 & v_2 & u_2^2 & u_2 v_2 & v_2^2 \\ u_3 & v_3 & u_3^2 & u_3 v_1 & v_3^2 \\ u_4 & v_4 & u_4^2 & u_4 v_4 & v_4^2 \\ \vdots & \vdots & \vdots & \vdots & \vdots \\ u_J & v_J & u_J^2 & u_J v_J & v_J^2 \end{bmatrix} \begin{bmatrix} \phi_x \\ \phi_y \\ \tfrac{1}{2}\phi_{xx} \\ \phi_{xy} \\ \tfrac{1}{2}\phi_{yy} \end{bmatrix} = \begin{bmatrix} \phi_1-\phi_i \\ \phi_2-\phi_i \\ \phi_3-\phi_i \\ \phi_4-\phi_i \\ \vdots \\ \phi_J-\phi_i \end{bmatrix},
\qquad (13)
$$

where

$$
u_j = u_j - u_i, \qquad v_j = v_j - v_i.
\qquad (14)
$$

Applying the least square method (LSM), we obtain

$$
\begin{bmatrix} u_1 & u_2 & u_3 & u_4 & \cdots & u_J \\ v_1 & v_2 & v_3 & v_4 & \cdots & v_J \\ u_1^2 & u_2^2 & u_3^2 & u_4^2 & \cdots & u_J^2 \\ u_1 v_1 & u_2 v_2 & u_3 v_3 & u_4 v_4 & \cdots & u_J v_J \\ v_1^2 & v_2^2 & v_3^2 & v_4^2 & \cdots & v_J^2 \end{bmatrix} \begin{bmatrix} W_1 & 0 & 0 & 0 & \cdots & 0 \\ 0 & W_2 & 0 & 0 & \cdots & 0 \\ 0 & 0 & W_3 & 0 & \cdots & 0 \\ 0 & 0 & 0 & W_4 & \cdots & 0 \\ \vdots & \vdots & \vdots & \vdots & \cdots & 0 \\ 0 & 0 & 0 & 0 & \cdots & W_J \end{bmatrix} \begin{bmatrix} u_1 & v_1 & u_1^2 & u_1 v_1 & v_1^2 \\ u_2 & v_2 & u_2^2 & u_2 v_2 & v_2^2 \\ u_3 & v_3 & u_3^2 & u_3 v_1 & v_3^2 \\ u_4 & v_4 & u_4^2 & u_4 v_4 & v_4^2 \\ \vdots & \vdots & \vdots & \vdots & \vdots \\ u_J & v_J & u_J^2 & u_J v_J & v_J^2 \end{bmatrix} \begin{bmatrix} \phi_x \\ \phi_y \\ \tfrac{1}{2}\phi_{xx} \\ \phi_{xy} \\ \tfrac{1}{2}\phi_{yy} \end{bmatrix}
$$

$$
= \begin{bmatrix} u_1 & u_2 & u_3 & u_4 & \cdots & u_J \\ v_1 & v_2 & v_3 & v_4 & \cdots & v_J \\ u_1^2 & u_2^2 & u_3^2 & u_4^2 & \cdots & u_J^2 \\ u_1 v_1 & u_2 v_2 & u_3 v_3 & u_4 v_4 & \cdots & u_J v_J \\ v_1^2 & v_2^2 & v_3^2 & v_4^2 & \cdots & v_J^2 \end{bmatrix} \begin{bmatrix} W_1 & 0 & 0 & 0 & \cdots & 0 \\ 0 & W_2 & 0 & 0 & \cdots & 0 \\ 0 & 0 & W_3 & 0 & \cdots & 0 \\ 0 & 0 & 0 & W_4 & \cdots & 0 \\ \vdots & \vdots & \vdots & \vdots & \cdots & 0 \\ 0 & 0 & 0 & 0 & \cdots & W_J \end{bmatrix} \begin{bmatrix} \phi_1-\phi_i \\ \phi_2-\phi_i \\ \phi_3-\phi_i \\ \phi_4-\phi_i \\ \vdots \\ \phi_J-\phi_i \end{bmatrix},
\qquad (15)
$$

where

$$
W_j = \begin{cases} \left( \dfrac{|\mathbf{r}_e - \mathbf{r}_i|}{|\mathbf{r}_j - \mathbf{r}_i|} - 1 \right) & |\mathbf{r}_j - \mathbf{r}_i| < |\mathbf{r}_e - \mathbf{r}_i| \\ 0 & \text{otherwise} \end{cases}.
\qquad (16)
$$

## 4.1. Verification in 2D using power function

In this example, the effects of weight are also studied. We use a following power function:

$$f(x,y) = x^4 + y^4 + x^3 y^3 \quad \text{in} \quad -2 < x < 2, \quad -2 < y < 2 \tag{17a}$$

$$f_x(x,y) = 4x^3 + 3x^2 y^3, \quad f_y(x,y) = 4y^3 + 3x^3 y^2, \tag{17b}$$

$$f_{xx}(x,y) = 12x^2 + 6xy^3, \quad f_{xy}(x,y) = 9x^2 y^2, \quad f_{yy}(x,y) = 12y^2 + 6x^3 y. \tag{17c}$$

Function $f(x,y)$ is calculated on $50 \times 50$ regular mesh points and is plotted in Fig. 1. The irregular mesh points $(x_i, y_j)$, $0 \le i, j \le 50$ are obtained by adding random number to the regular mesh points $(x_{i\_regular}, y_{j\_regular})$, $0 \le i, j \le 50$:

$$\begin{aligned} x_i &= x_{i\_regular} + \Delta r \times (\text{uniform random number in } (-1, +1)) \\ y_j &= y_{j\_regular} + \Delta r \times (\text{uniform random number in } (-1, +1)) \end{aligned} \tag{18}$$

In the following examples, the scale of irregularity $\Delta r$ and the radius $r$ defining the neighboring points are $0.25\Delta x$ and $2.5\Delta x$, respectively. The irregular mesh points are shown in Fig. 2.

Figures 3, 4 and 5 show the comparisons among the analytical result, FD (Finite Difference) result, DDIN (Discrete Differential operators on Irregular Nodes) result without weight and DDINW result with weight of $f_x(x,y)$, $f_{xx}(x,y)$ and $f_{xy}(x,y)$, respectively. The effects of weight are shown in Figs. 3, 4 and 5. The effects are not big in this example.

The rms (root mean square) errors are given in Fig. 6. In Fig. 6, Line $x^1$ is proportional to $\Delta x$, and line $x^2$ is proportional to $(\Delta x)^2$. The definition of rms is given by, for example

$$\text{rms of } f_x = \sqrt{\frac{1}{\text{No. of total points}} \sum_{j=1}^{J} (f_{x\_numerical}(x_i, y_i) - f_{x\_exact}(x_i, y_i))^2} \tag{19}$$

In Fig. 6, the effect of data spacing $\Delta x$ on the root mean square (RMS) errors of $f_x(x,y)$, $f_{xx}(x,y)$ and $f_{xy}(x,y)$ are shown. The error estimation given by Eq. (11) is verified by numerical calculations, since the rms errors of $f_x(x,y)$ is proportional to $(\Delta x)^2$ and those of $f_{xx}(x,y)$ and $f_{xy}(x,y)$ are proportional to $\Delta x$

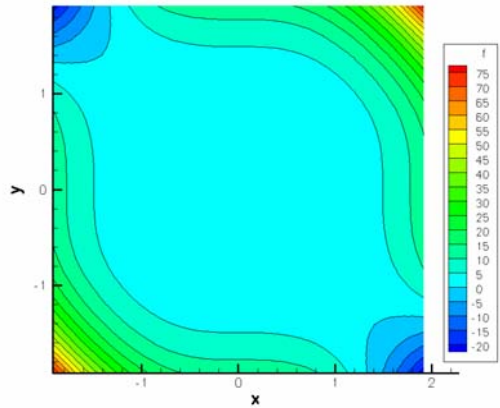

Fig. 1. Function $f(x,y)$.

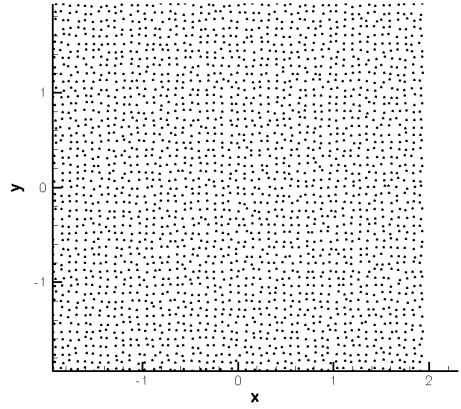

Fig. 2. Distribution of irregular node points ($\Delta r = 0.25\Delta x$).

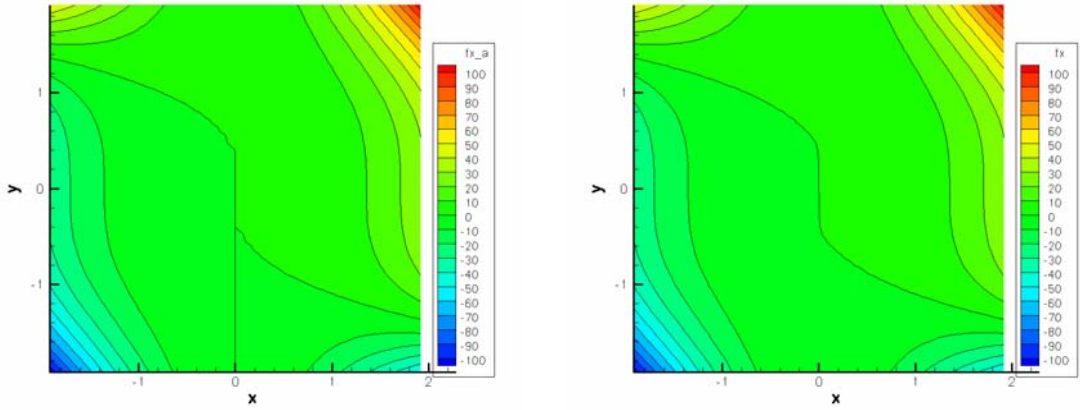

(a) $f_x(x,y)$: analytical    (b) $f_x(x,y)$: FD

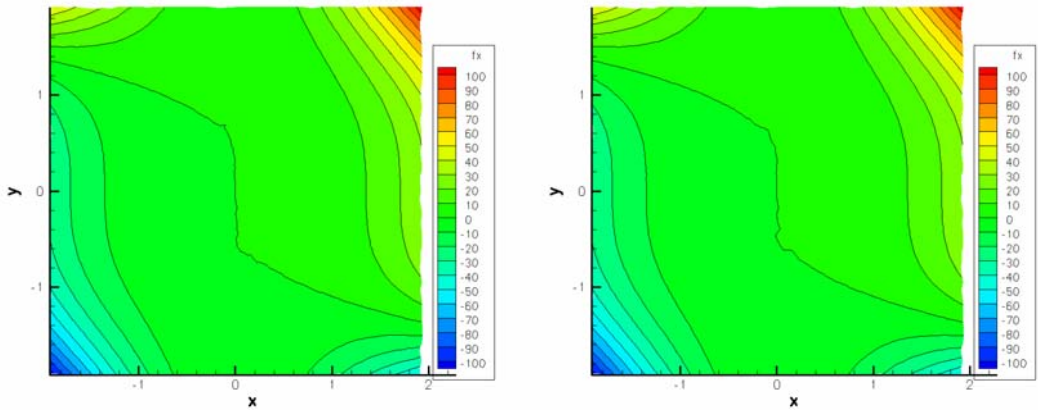

(c) $f_x(x,y)$: DDIN ($\Delta r = 0.25\Delta x, r = 2.5\Delta x$)    (d) $f_x(x,y)$: DDINW ($\Delta r = 0.25\Delta x, r = 2.5\Delta x$)

Fig. 3. Comparison of $f_x(x,y)$ among analytical, finite difference, DDIN and DDINW.

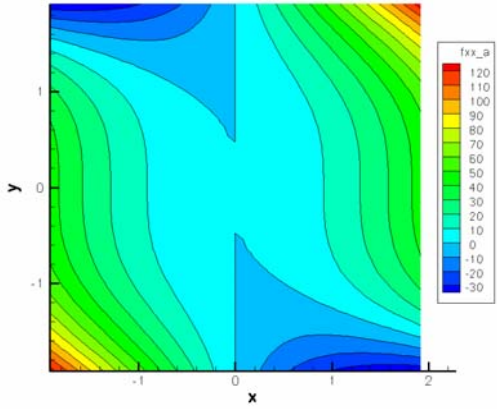
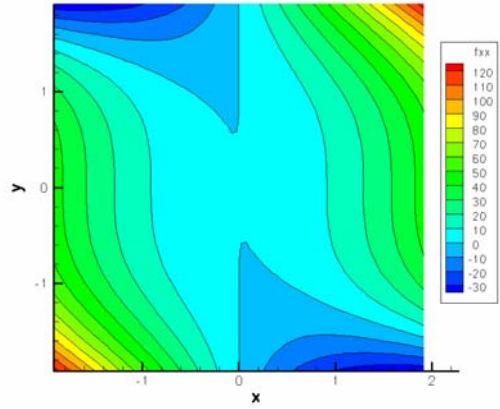

(a) $f_{xx}(x,y)$: analytical  (b) $f_{xx}(x,y)$: FD

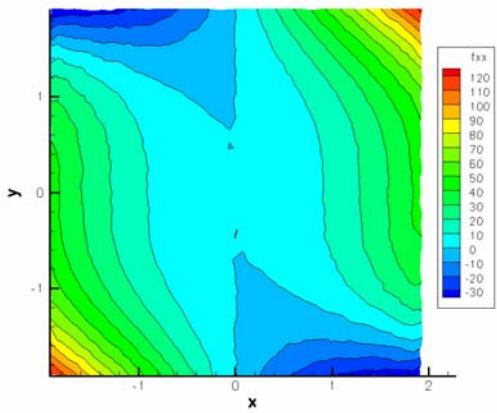
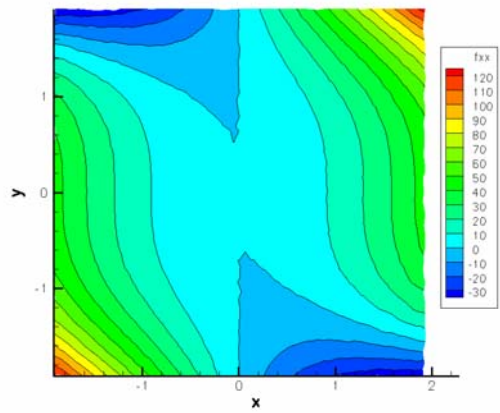

(c) $f_{xx}(x,y)$: DDIN ($\Delta r = 0.25\Delta x$, $r = 2.5\Delta x$)  (d) $f_{xx}(x,y)$: DDINW ($\Delta r = 0.25\Delta x$, $r = 2.5\Delta x$)

Fig. 4. Comparison of $f_{xx}(x,y)$ among analytical, finite difference, DDIN and DDINW.

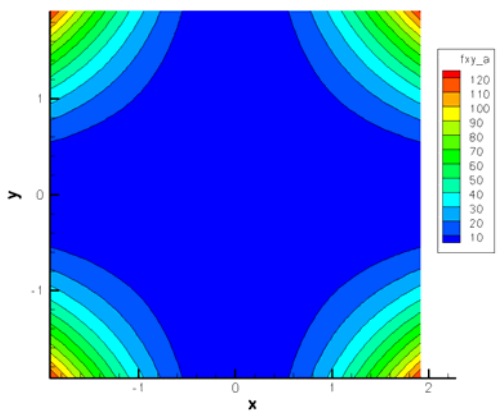
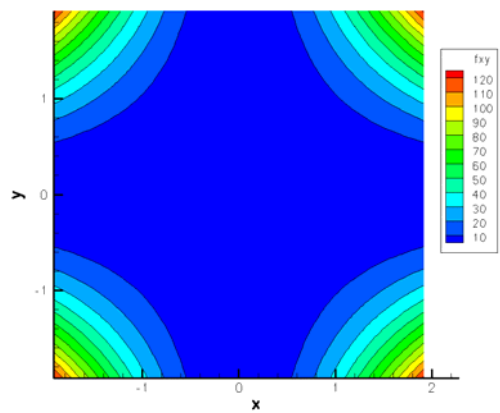

(a) $f_{xy}(x,y)$: analytical  (b) $f_{xy}(x,y)$: FD

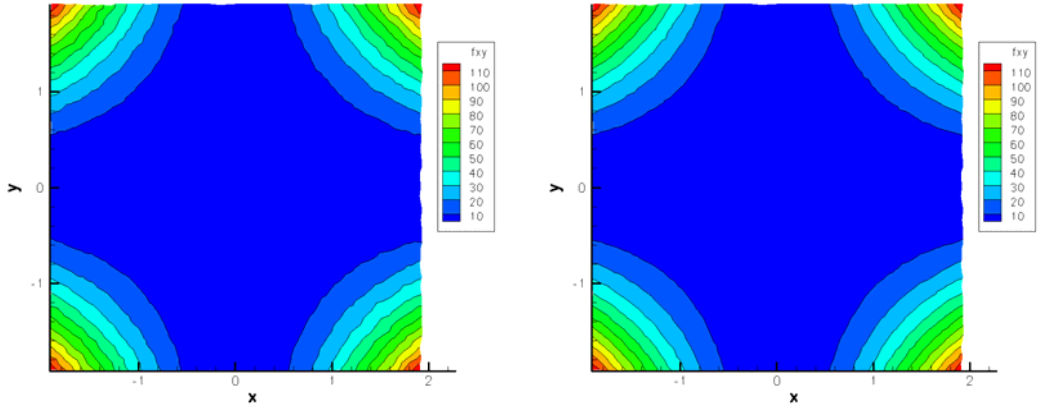

(c) $f_{xy}(x,y)$: DDIN ($\Delta r = 0.25\Delta x, r = 2.5\Delta x$)   (d) $f_{xy}(x,y)$: DDINW ($\Delta r = 0.25\Delta x, r = 2.5\Delta x$)

Fig. 5. Comparison of $f_{xy}(x,y)$ among analytical, finite difference, DDIN and DDINW.

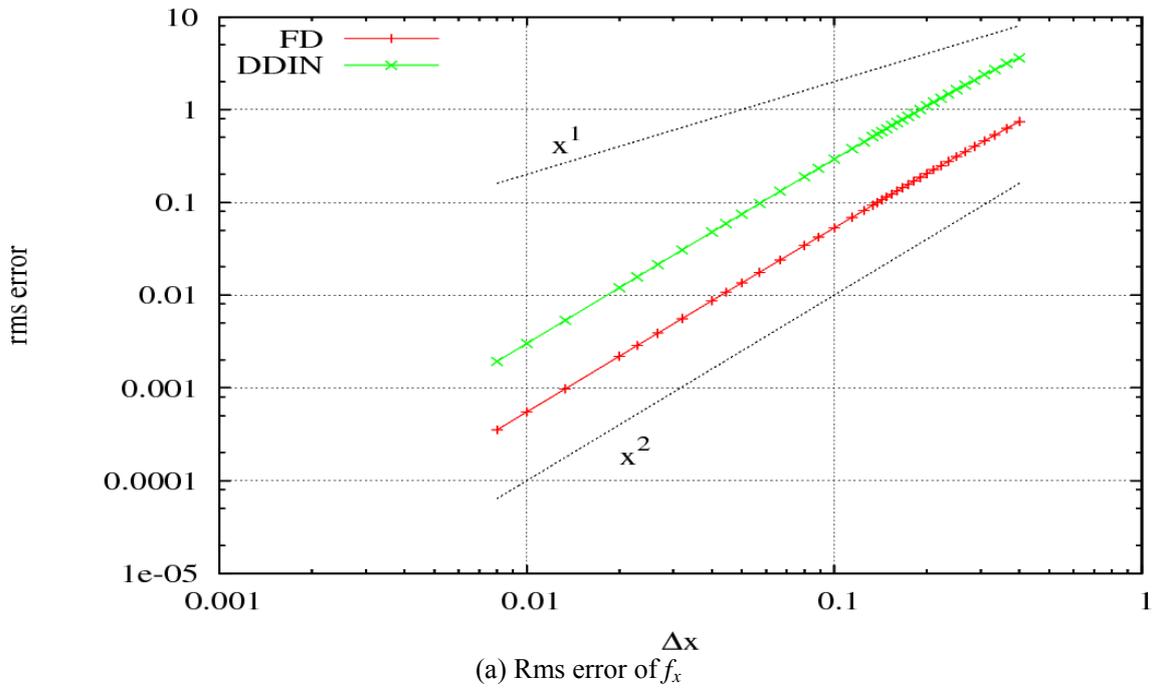

(a) Rms error of $f_x$

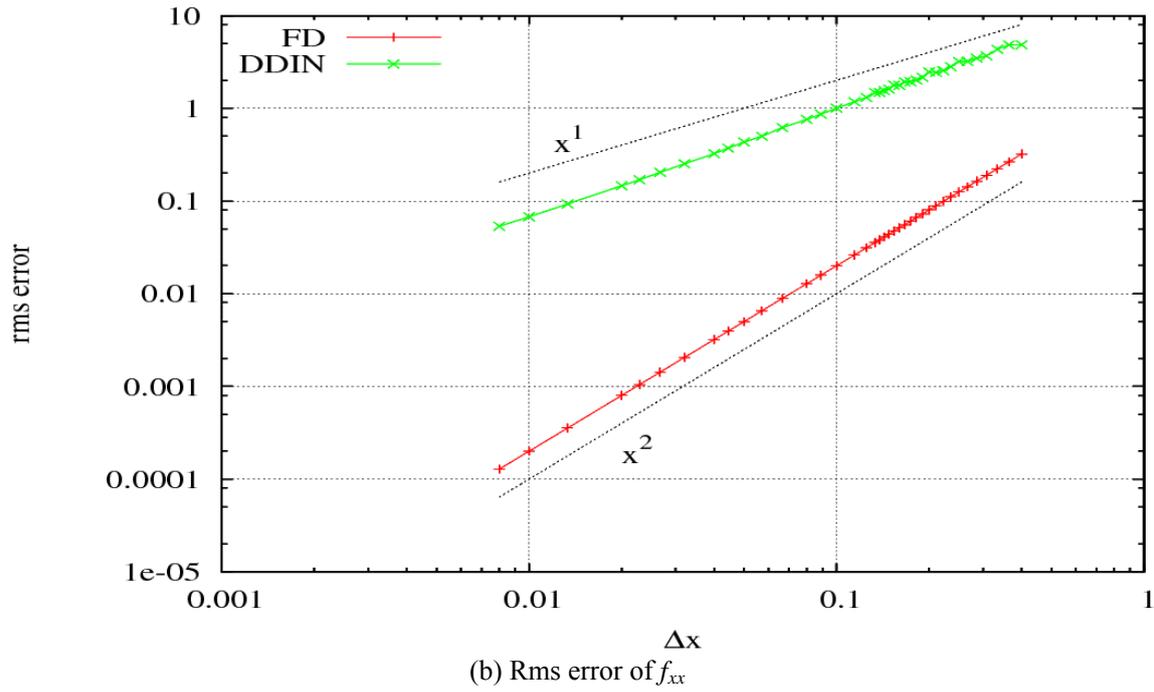
(b) Rms error of $f_{xx}$

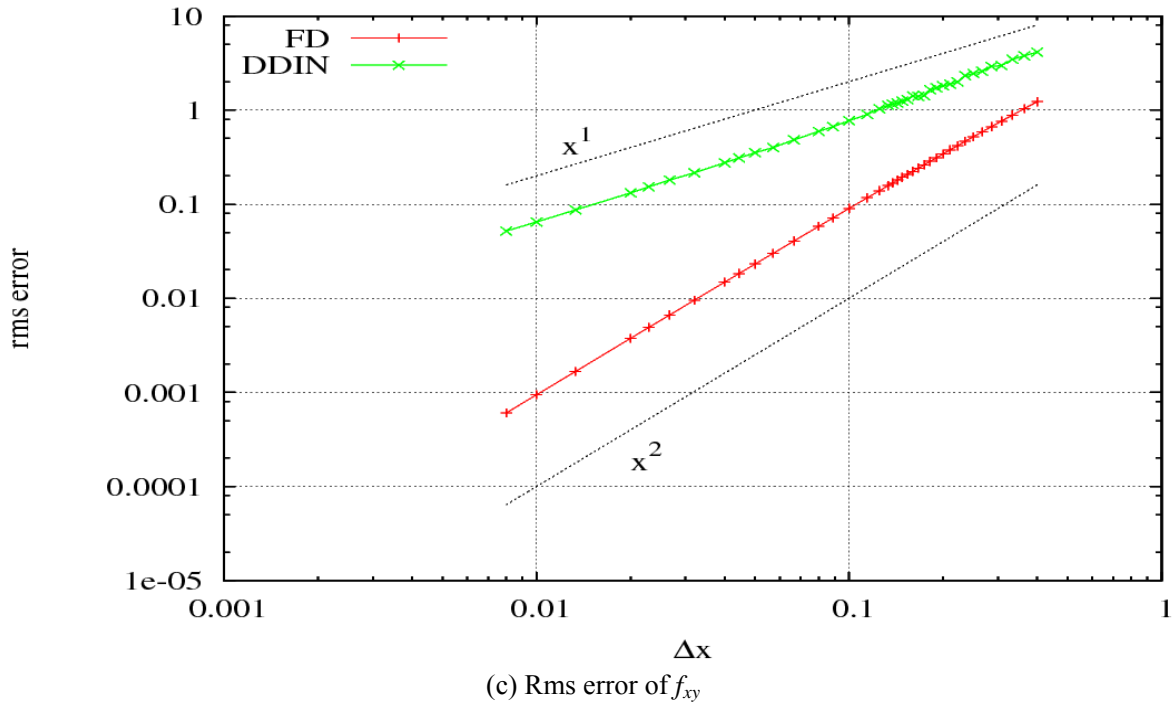
(c) Rms error of $f_{xy}$

Fig. 6. Comparison of $f_{xy}(x, y)$ between FD and DDIN.

**4.2. Verification in 2D using sinusoidal function with weight**

In this example, the effects of weight are also studied. We use a sinusoidal function $f(x, y)$:

$$f(x, y) = \sin x \cos 2y \text{ in } 0 \leq x \leq 2\pi, 0 \leq y \leq 2\pi, \quad (20a)$$

$$f_x(x, y) = \cos x \cos 2y, \quad f_y(x, y) = -2 \sin x \sin 2y, \quad (20b)$$

$$f_{xx}(x, y) = -\sin x \cos 2y, \quad f_{xy}(x, y) = -2 \cos x \sin 2y, \quad f_{yy}(x, y) = -4 \sin x \cos 2y. \quad (20c)$$

Function $f(x, y)$ is calculated on $50 \times 50$ regular mesh points and is plotted in Fig. 7. The irregular mesh points $(x_i, y_j)$, $0 \leq i, j \leq 50$ are obtained by adding random number to the regular mesh points

$\left(x_{i\_regular}, y_{j\_regular}\right)$, $0 \leq i, j \leq 50$:

$$x_i = x_{i\_regular} + \Delta r \times (uniform\ random\ number\ in\ (-1,+1))$$
$$y_j = y_{j\_regular} + \Delta r \times (uniform\ random\ number\ in\ (-1,+1))$$ . (21)

In the following examples, the scale of irregularity $\Delta r$ and the radius $r$ defining the neighboring points are $0.25\Delta x$ and $3\Delta x$, respectively. The irregular mesh points are shown in Fig. 8.

Figures 9 through 13 show the comparisons among the analytical result, FD (Finite Difference) result, DDIN (Discrete Differential operators on Irregular Nodes) result without weight and DDINW result with weight of $f_x(x,y)$, $f_y(x,y)$, $f_{xx}(x,y)$, $f_{xy}(x,y)$ and $f_{yy}(x,y)$, respectively. The effects of weight are clearly shown in Figs. 11, 12 and 13.

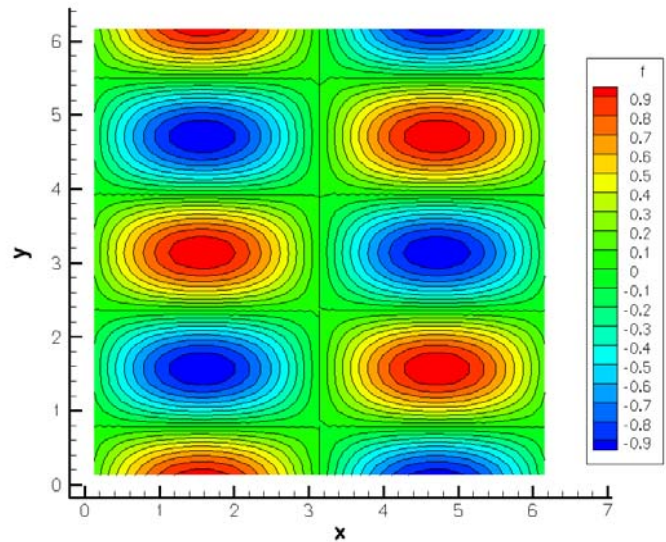

Fig. 7. Function f(x, y).

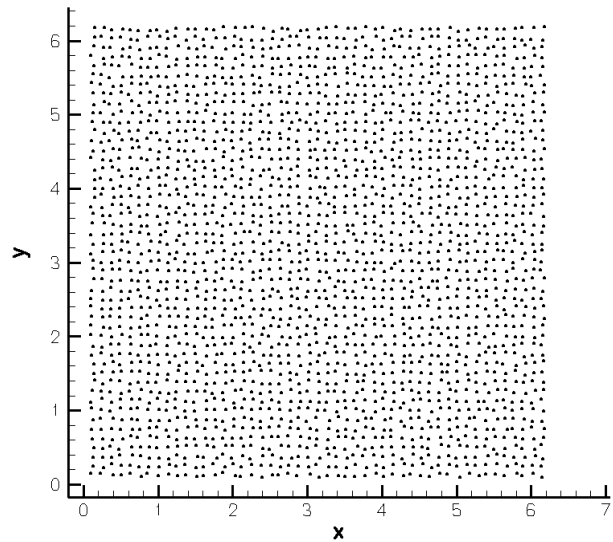

Fig. 8. Distribution of points( $\Delta r = 0.25\Delta x$ )

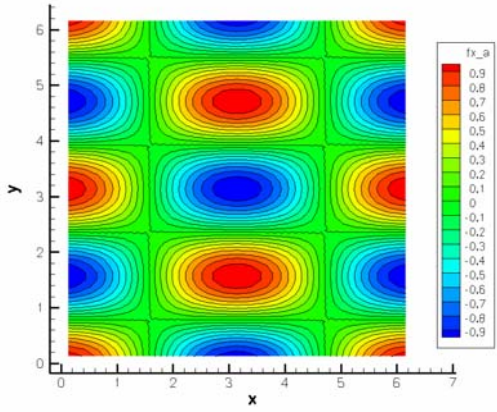
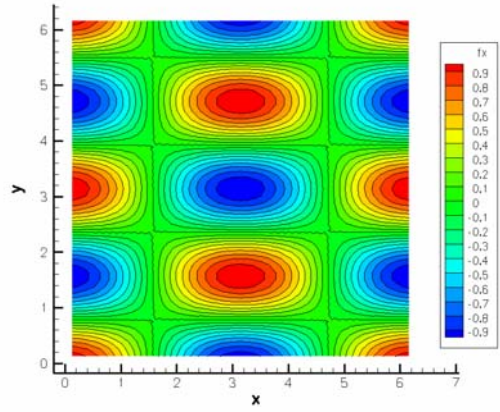

$f_x(x, y)$: analytical          $f_x(x, y)$: FD

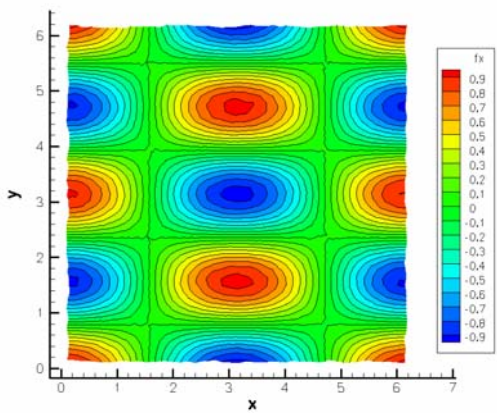
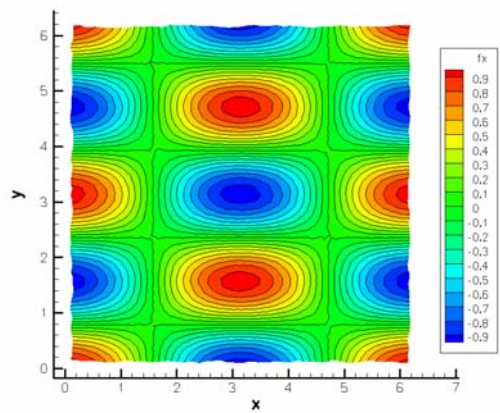

$f_x(x, y)$: DDIN ( $\Delta r = 0.25\Delta x$, $r = 3\Delta x$ )     $f_x(x, y)$: DDINW ( $\Delta r = 0.25\Delta x$, $r = 3\Delta x$ )

Fig. 9. Comparison of $f_x(x, y)$

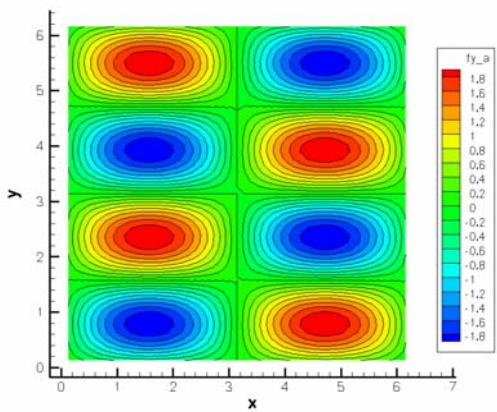
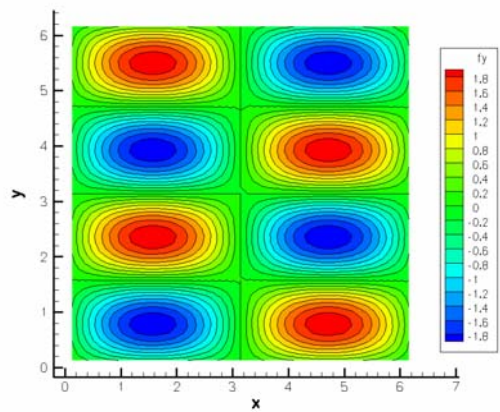

$f_y(x, y)$: analytical          $f_y(x, y)$: FD

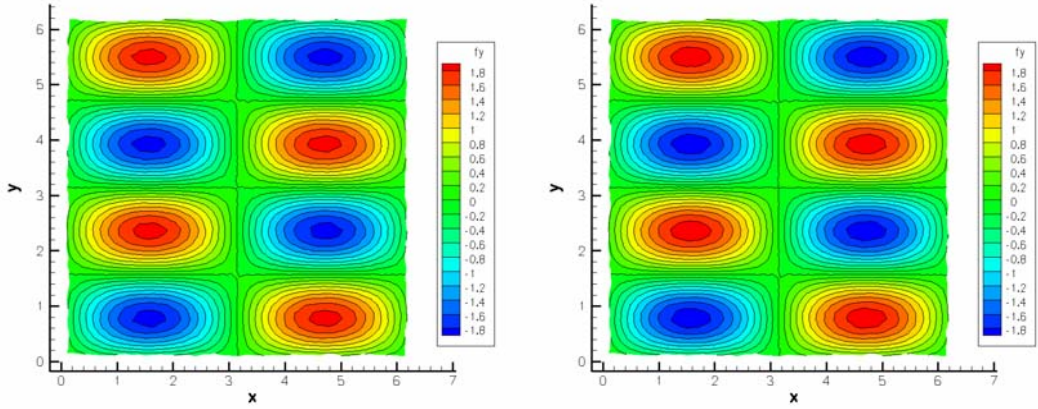

$f_y(x, y)$: DDIN  $f_y(x, y)$: DDINW

($\Delta r = 0.25\Delta x$, $r = 3\Delta x$)

Fig. 10. Comparison of $f_y(x, y)$

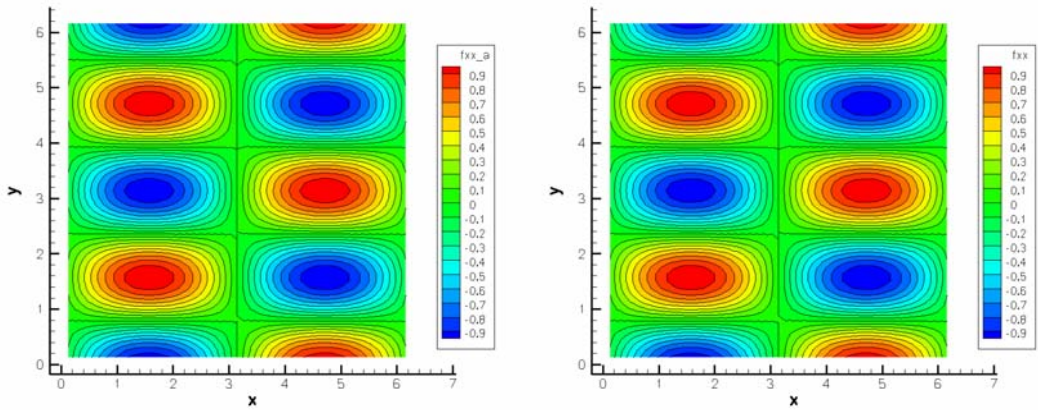

$f_{xx}(x, y)$: analytical  $f_{xx}(x, y)$: FD

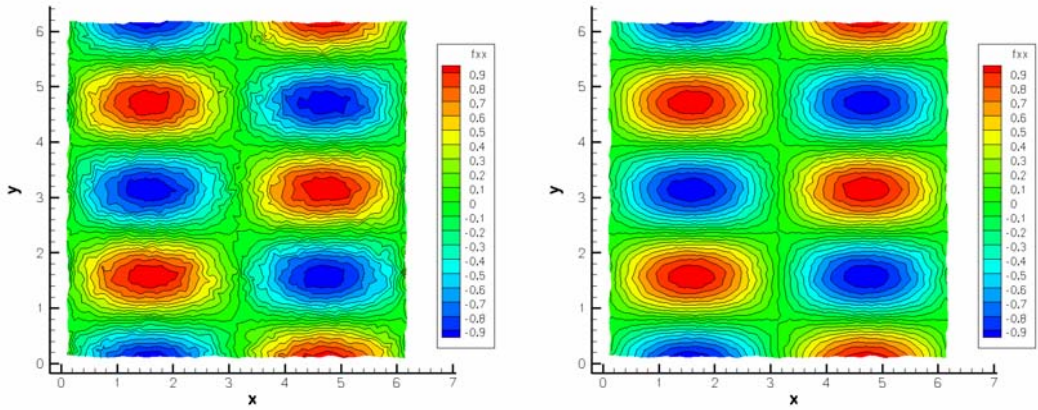

$f_{xx}(x, y)$: DDIN  $f_{xx}(x, y)$: DDINW

($\Delta r = 0.25\Delta x$, $r = 3\Delta x$)

Fig. 11. Comparison of $f_{xx}(x, y)$

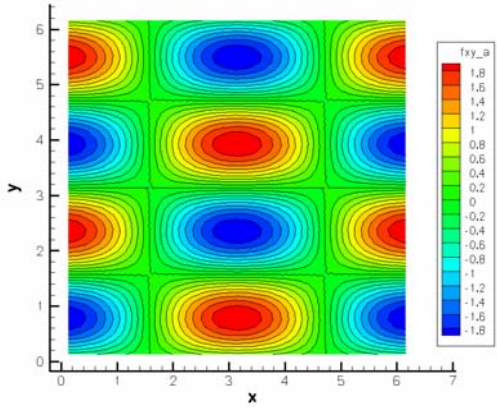
$f_{xy}(x, y)$: analytical

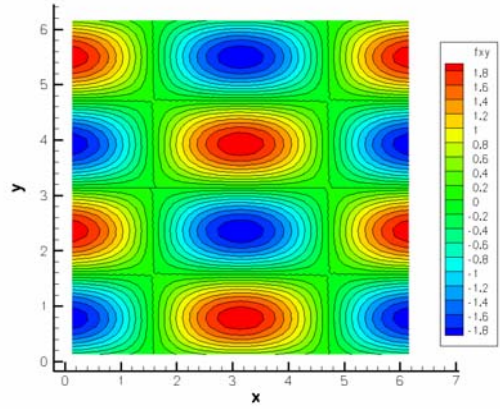
$f_{xy}(x, y)$: FD

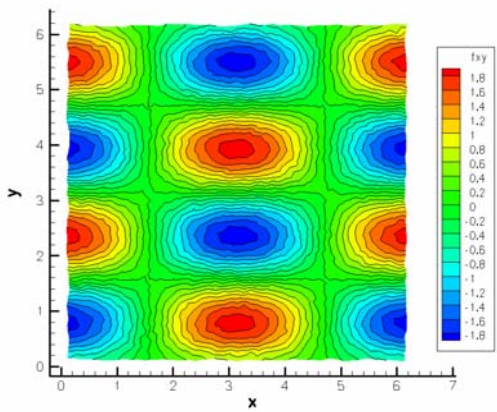
$f_{xy}(x, y)$: DDIN

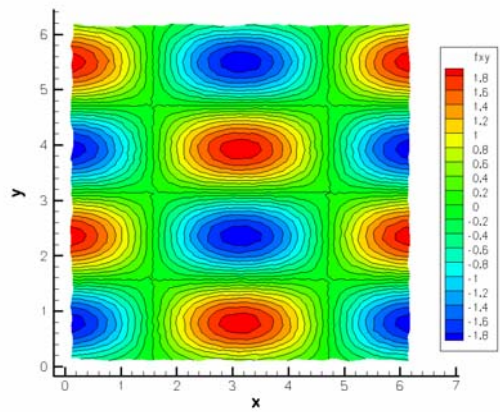
$f_{xy}(x, y)$: DDINW

($\Delta r = 0.25\Delta x$, $r = 3\Delta x$)
Fig. 12. Comparison of $f_{xy}(x, y)$

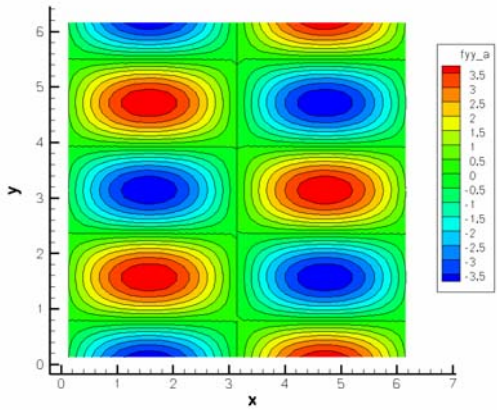
$f_{yy}(x, y)$: analytical

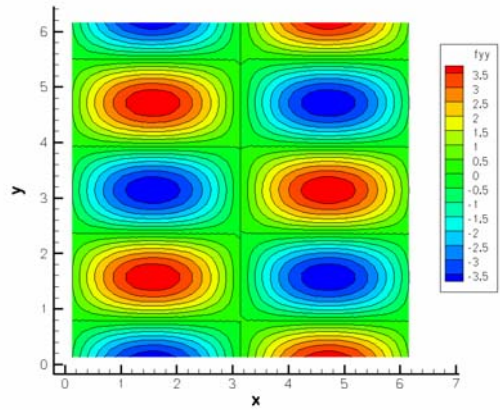
$f_{yy}(x, y)$: FD

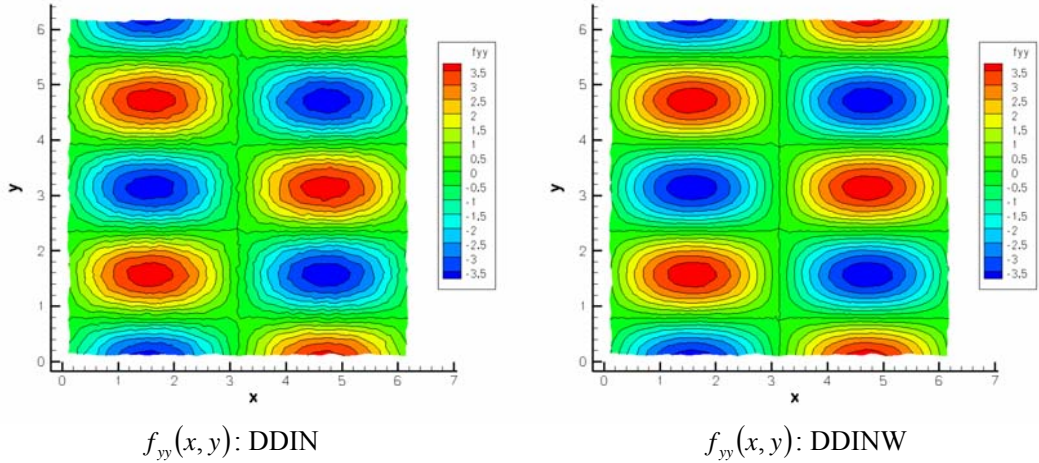

$f_{yy}(x,y)$: DDIN    $f_{yy}(x,y)$: DDINW

($\Delta r = 0.25\Delta x$, $r = 3\Delta x$)

Fig. 13. Comparison of $f_{yy}(x,y)$

## 5. Conclusions

We have discussed the error estimation of difference operators in irregular nodes theoretically. We verified our theory by numerical calculations.

If we compare Eqs. (1) and (11), the accuracies of the first order difference operators are same unexpectedly. If we compare Eqs. (4) and (11), the accuracy becomes lower in irregular nodes than in regular nodes. The symmetry of the node allocation in regular nodes increases the accuracy.

We also showed the interesting effects of the weight on the second order difference.